\documentclass[amssymb,amstex,amsmath,reqno, 11pt]{amsart}

\usepackage{graphicx}
\usepackage{amsmath}
\usepackage{amssymb}

\newtheorem{thm}{Theorem}
\newenvironment{pf}{\noindent \rm Proof.}

\begin{document}

\title[new minimal surfaces in $\mathbb S^3$]
{NEW MINIMAL SURFACES IN $\mathbb S^3$ desingularizing the Clifford tori}
\author[J. CHOE and M. Soret]{JAIGYOUNG CHOE and MARC SORET}
\date{March 9, 2013}
\thanks{J.C. supported in part by NRF, 2011-0030044, SRC-GAIA}
\address{Jaigyoung Choe: Korea Institute for Advanced Study, Seoul, 130-722, Korea \\
Marc Soret: D\'{e}partement de Math\'{e}matiques, Universit\'{e} F. Rabelais, 37200 Tours,
France}
\email{choe@kias.re.kr, ma.soret@noos.fr}

\begin{abstract}
For each integer $m\geq2$ and $\ell\geq1$ we construct a pair of compact embedded minimal
surfaces of genus $1+4m(m-1)\ell$. These surfaces desingularize the $m$
Clifford tori meeting each other along a great circle at the angle of
$\pi/m$. They are invariant under a finite group of screw motions and have no reflection symmetry across a great sphere.
\end{abstract}

\maketitle

One can construct a complete minimal surface in $\mathbb R^3$ by suitably choosing a
holomorphic 1-form $f(z)dz$ and a meromorphic function $g(z)$ for the Weierstrass
representation formula. But there is no such an efficient tool in $\mathbb S^3$. This is why
not many minimal surfaces are known to exist in $\mathbb S^3$. So far only three types of
minimal surfaces have been constructed and their methods of construction are all different.
In 1970 Lawson \cite{L} constructed infinitely many compact minimal surfaces in $\mathbb
S^3$; in 1988 Karcher-Pinkall-Sterling \cite{KPS} found nine new compact embedded minimal
surfaces; in 2010 Kapouleas-Yang \cite{KY} obtained new minimal surfaces by doubling the
Clifford torus.

Lawson starts from a piecewise geodesic Jordan curve $\Gamma$, finds a minimal disk $D$
spanning $\Gamma$, and extends $D$ across $\Gamma$ by $180^\circ$-rotations to obtain a compact
immersed minimal surface. Lawson's Jordan curve $\Gamma$ consists of 4 geodesic segments and
is a subset of the 1-skeleton of a tetrahedron in $\mathbb S^3$. This tetrahedron is a
fundamental piece of a tessellation of $\mathbb S^3$.

On the other hand, Karcher-Pinkall-Sterling start from a tetrahedron $T$ which gives rise to
a different type of tessellation of $\mathbb S^3$. Then they find a minimal disk $D$ in $T$
which is perpendicular to $\partial T$ along $\partial D$, and extend $D$ by the reflections across $\partial T$ to obtain a compact embedded minimal surface.

Kapouleas-Yang's minimal surfaces resemble two parallel copies of the Clifford torus, joined
by $m^2$ small catenoidal bridges for sufficiently large $m$ symmetrically arranged along a
square lattice of points on the torus.

In this paper we construct infinitely many compact embedded minimal surfaces by desingularizing $m$ Clifford tori which meet each other along a great circle at the angle of $\pi/m$. Our desingularization does not employ the gluing method, instead we use a tessellation of $\mathbb S^3$ by $16m^2\ell$ ($m\geq2, \ell\geq1$) pentahedra and apply Lawson's method for the Jordan curve of 6 geodesic segments which is a subset of the 1-skeleton of a pentahedron. The resulting compact embedded minimal surface has genus $1+4m(m-1)\ell$ (Theorem 1).

Given a great circle $C_1$ in $\mathbb S^3$, there is the polar great circle $C_2$ of $C_1$, that is, dist$(p,q)=\pi/2$ for any $p\in C_1$ and $q\in C_2$. $C_1$ and $C_2$ are linked in $\mathbb S^3$. If $m$ Clifford tori meet each other along $C_1$, then they intersect each other along $C_2$ as well. Therefore once $m$ Clifford tori are desingularized along $C_1$, there are two ways of desingularizing the tori along $C_2$. Thus we obtain the second type(even) of minimal surfaces desingularizing $m$ Clifford tori for each genus $1+4m(m-1)\ell,\ell\geq2$ (Theorem 2).

All the embedded minimal surfaces constructed by Lawson, Karcher-Pinkall-Sterling, Kapouleas-Yang satisfy the reflection symmetry, i.e., they are invariant under a reflection across a great sphere in $\mathbb S^3$. But our new minimal surfaces have no reflection symmetry.

\section{Clifford torus}
The Clifford torus $T$ is the building block of our new minimal surfaces. So we start by investigating its two characteristic properties: it has the {\it equidistance property} and is {\it doubly ruled}. Define
$$T=\mathbb S^1(1/\sqrt{2})\times\mathbb S^1(1/\sqrt{2})=\{(x_1,x_2,x_3,x_4)\in\mathbb R^4:x_1^2+x_2^2=x_3^2+x_4^2=1/2\}.$$
Let $C_{12}, C_{34}$ be the two linked great circles in $\mathbb S^3$ defined by $C_{12}=\{(x_1,x_2,0,0):x_1^2+x_2^2=1\}$, $C_{34}=\{(0,0,x_3,x_4):x_3^2+x_4^2=1\}$.
Throughout this paper ``dist" denotes the distance in $\mathbb S^3$. Then
$${\rm dist}(p,q)=\pi/2,\,\,\forall p\in C_{12},\,\,\forall q\in C_{34},$$
and one gets the equidistance property:
$${\rm dist}(T,C_{12})={\rm dist}(T,C_{34})=\pi/4.$$
Also it is easy to see that
$$\overline{pq}\,\perp\, T,\,\,\forall p\in C_{12},\,\,\forall q\in C_{34}.$$

Let $\gamma_1=\{(x_1,x_2,1/\sqrt{2},0):x_1^2+x_2^2=1/2\},\,\gamma_2=\{(1/\sqrt{2},0,x_3,x_4):x_3^2+
x_4^2=1/2\}$. Cutting out $\gamma_1$ and $\gamma_2$ from $T$, one can obtain a flat square $Q\subset T$. Then one can
 consider two 1-parameter families of lines on $T$ which are parallel to the two diagonals of the square $Q$. These lines of $T$ are in fact the great circles of $\mathbb S^3$. For this reason $T$ is called doubly ruled. Let's see why these lines are great circles.
$$x_1^2+x_2^2=x_3^2+x_4^2\,\,{\rm becomes}\,\,(x_1+x_3)(x_1-x_3)=(x_4+x_2)(x_4-x_2).$$
Hence if we rotate $x_1x_3$-plane and $x_2x_4$-plane by $\pi/4$ and by $-\pi/4$, respectively, and use $x_1,x_2,x_3,x_4$ again for the new coordinates, then we get
$$x_1x_3=x_2x_4.$$
Hence $T$ can be represented by the coordinate map $\Psi:[0,2\pi)\times[0,2\pi)\rightarrow\mathbb S^3$,
$$\Psi(x,y)=(\cos x\sin y, \cos x\cos y,\sin x\cos y,\sin x\sin y).$$
Here we claim that $T$ is ruled by the two families of great circles $\{x={\rm const}\}$ and $\{y={\rm const}\}$.

Let $\rho_{ij}^t$ be the counterclockwise rotation of $\mathbb S^3$ by the angle $t$ along the $x_ix_j$-plane and define
$$\Phi_{ijkl}^{t}=\rho_{ij}^t\circ\rho_{kl}^{t},$$
where $\{i,j,k,l\}=\{1,2,3,4\}$ as a set. We will call $\Phi_{ijkl}^{t}$ a {\it screw motion} because it can be viewed as the composition of a rotation and a translation, $\rho_{kl}^{t}$ being the translation along the great circle $x_k^2+x_l^2=1$. Note that
\begin{eqnarray*}
\Psi(x,y)&=&\cos x(\sin y,\cos y,0,0)+\sin x(0,0,\cos y,\sin y)\\
&=&\cos y(0,\cos x,\sin x,0)+\sin y(\cos x,0,0,\sin x).
\end{eqnarray*}
Hence $T$ is foliated by the great circles $\{\Phi_{1423}^{t}(C_{21})\}$ which are $\{x={\rm const}\}$ and by the great circles $\{\Phi_{2134}^{t}(C_{23})\}=\{y={\rm const}\}$. Here $C_{21}$ is the great circle $C_{12}$ with the opposite orientation and $C_{23}=\{(0,x_2,x_3,0):x_2^2+x_3^2=1\}$.

These two families of great circles are orthogonal to each other. The orthogonality can be observed more easily on the fundamental piece $\hat{T}$ of the Clifford torus $T$ as in Figure 1. $T$ consists of eight congruent pieces of $\hat{T}$ and $\hat{T}$ is Morrey's solution to the Plateau problem for a geodesic polygon $\Gamma=\hat{C}_{12}\cup\hat{C}_{23}\cup\hat{C}_{34}\cup\hat{C}_{41}$ as used by Lawson in \cite{L}. $\hat{C}_{ij}$ is a subarc of length $\pi/2$ of $C_{ij}$ such that $\hat{C}_{12}$ is from $(1,0,0,0)$ to $(0,1,0,0)$, $\hat{C}_{23}$ from $(0,1,0,0)$ to $(0,0,1,0)$, $\hat{C}_{34}$ from $(0,0,1,0)$ to $(0,0,0,1)$, and $\hat{C}_{41}$ from $(0,0,0,1)$ to $(1,0,0,0)$.

\begin{center}
\includegraphics[width=3in]{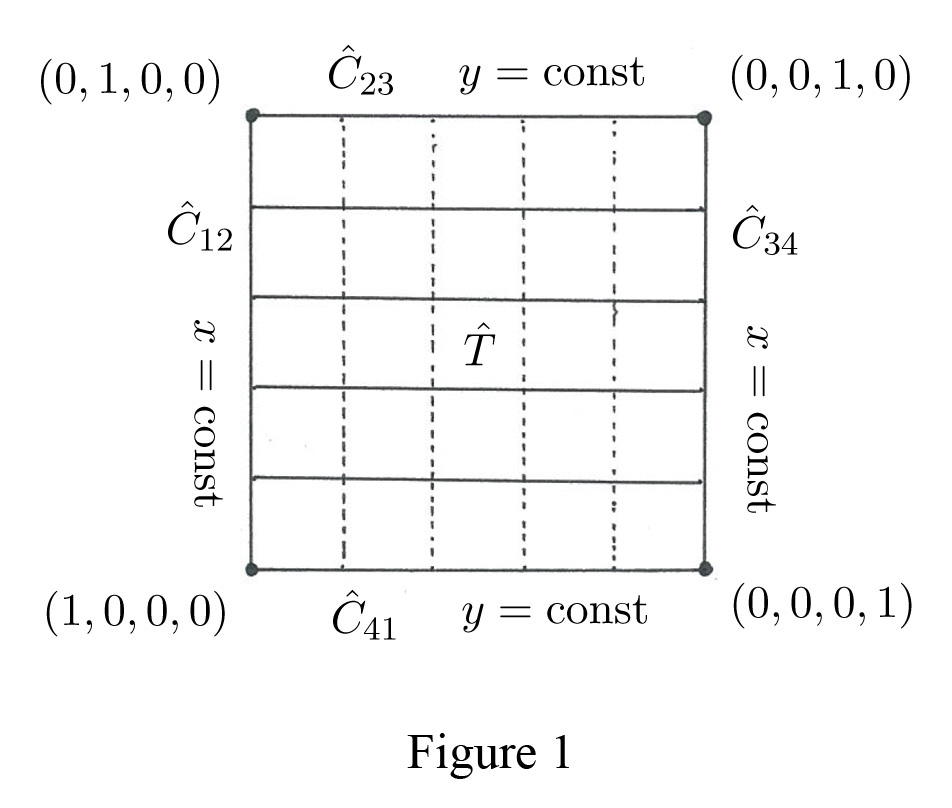}\\
\end{center}

Finally it is not difficult to see that $T$ is the equidistance set from the two great circles $\Phi_{1342}^{\pi/4}(C_{12})$ and $\Phi_{1342}^{3\pi/4}(C_{12})$ (remember that the original $x_1x_3$-plane and $x_2x_4$-plane have been rotated by $\pi/4$ and $-\pi/4$, respectively). Also it should be mentioned that $T$ is invariant under the screw motions $\Phi_{1234}^{t}$ and $\Phi_{1423}^{t}$ for any $t$. And if a great circle of $\mathbb S^3$ lies in a Clifford torus, so does its polar circle.

\section{Odd surfaces}
Given two orthogonal planes in $\mathbb R^3$, the minimal surface that desingularizes them is Scherk's second surface. For two great spheres orthogonal to each other in $\mathbb S^3$, the minimal surfaces that desingularize them are Lawson's minimal surfaces $\xi_{m,k}$ of genus $mk$. Then, given two orthogonal Clifford tori in $\mathbb S^3$, is there a minimal surface that desingularizes them? We are motivated by this question and are led to the following.

\begin{thm} Let $T_1,\ldots,T_m\subset \mathbb S^3$ be the Clifford tori
    intersecting each other along a great circle $C_1$ at an angle of $\pi/m$. Then there
    exists a compact minimal surface $T_{m,k}^{o}$ desingularizing $T_1\cup\cdots\cup T_m$ for each
    $k=2m\ell$ with integer $\ell\geq1${\rm :}\\
    {\rm (i)} $T_{m,k}^{o}$ is embedded and has
    genus $1+2k(m-1)=1+4m(m-1)\ell${\rm ;}\\
    {\rm (ii)} $T_{m,k}^{o}$ is invariant under a finite group of screw motions{\rm ;}\\
    {\rm (iii)} $T_{m,k}^{o}$ has no reflection symmetry across a great sphere{\rm ;}\\
    {\rm (iv)} {\rm Area}$(T^o_{m,k})<2m\pi^2$.
\end{thm}
\begin{pf}
Let $C_2$ be the polar  great circle of $C_1$, that is, the set of all points of distance $\pi/2$ from $C_1$. Then $T_1\cap\cdots\cap T_m=C_1\cup C_2$. We claim that $T_1,\ldots,T_m$ also meet each other along $C_2$ at the angle of $\pi/m$. Introduce the coordinates $x_1,x_2,x_3,x_4$ of $\mathbb R^4\supset\mathbb S^3$ such that $C_1=C_{12}:x_1^2+x_2^2=1,\,C_2=C_{34}:x_3^2+x_4^2=1$. Then $T_1,\ldots,T_m$ are invariant under $\Phi_{1234}^{t}$ or $\Phi_{1243}^{t}$. Suppose without loss of generality that $\Phi_{1234}^{t}(T_1)=T_1$. Then $$\{T_1,\ldots,T_m\}=\{T_1,\rho_{34}^{\pi/m}(T_1),\rho_{34}^{2\pi/m}(T_1),\ldots,
\rho_{34}^{(m-1)\pi/m}(T_1)\}.$$ Since $T_1=\rho_{34}^{-t}\circ\rho_{12}^{-t}(T_1)$, it follows that $$\{T_1,\ldots,T_m\}=\{T_1,\rho_{12}^{-\pi/m}(T_1),\ldots,\rho_{12}^{-(m-1)\pi/m}(T_1)\},$$ hence the claim follows.

Let $T_0$ be the Clifford torus which is the equidistance set from $C_1$ and $C_2$. $T_0$ divides $\mathbb S^3$ into the two domains denoted $D_1$ and $D_2$ containing $C_1$ and $C_2$, respectively. Choose equally spaced points $p_1,\ldots,p_{2k}$ on $C_1$ such that ${\rm dist}(p_j,p_{j+1})=\pi/k$ and let $S_1^1,\ldots,S_{2k}^1$ be the great spheres such that $S_j^1$ contains $p_j$ and is perpendicular to $C_1$ at $p_j$. $T_1,\ldots,T_m$ and $S_1^1,\ldots,S_{2k}^1$ divide $D_1$ into congruent domains $\{U_j^i\}_{1\leq i\leq2m,\,1\leq j\leq2k}$. $\{U_j^i\}$ are numbered in such a way that $\cup_{i=1}^{2m}U_j^i$ is a component of $D_1\sim(S_j^1\cup S_{j+1}^1)$ and $\cup_{j=1}^{2k}U_j^i,\cup_{j=1}^{2k}U_j^{i+m}$ are components of $D_1\sim(T_i\cup T_{i+1})$, and $U_j^i,U_j^{i+m}$ are symmetric about $C_1$, that is, $U_j^{i+m}=\rho_{C_1}(U_j^i)$, $\rho_{C}$ denoting the $180^\circ$-rotation about the great circle $C$. These domains are in fact congruent pentahedra bounded by three Clifford tori and two great spheres as in Figure 2. Recall that the two great spheres are perpendicular to the base Clifford torus $T_0$. Each $\bar{U}_j^i\cap T_0$ is a parallelogram on $T_0$. Hence the tessellation of $D_1$ by the pentahedra $\{U_j^i\}$ gives rise to a tessellation of $T_0$ by the parallelograms $\{A_j^i\}_{1\leq i\leq2m,\,1\leq j\leq2k}$.

\begin{center}
\includegraphics[width=3in]{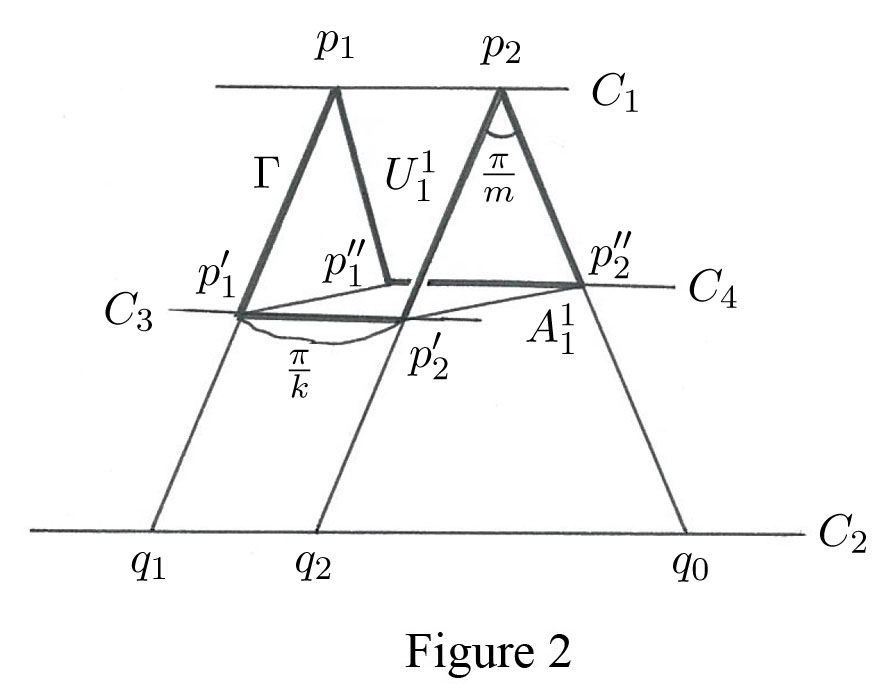}\\
\end{center}

Note that $U_1^1$ is bounded by $T_0,T_1,T_2$ and $S_1^1,S_2^1$, with $p_1'\in T_1,p_1''\in T_2$. Denote by $p_1',p_2',p_1'',p_2''$ the vertices of $A_1^1$ (also of $U_1^1$). Let $c(s)$ be the arclength parametrization of the geodesic $\overline{p_1p_1'}$ with $c(0)=p_1, c(\pi/4)=p_1'$. Then the angle between $T_1$ and $S_1^1$ at $c(s)$ equals $\pi/2-s$ because the tangent plane to $T_1$ at $c(s)$ is rotating around $\overline{p_1p_1'}$ under a screw motion as $s$ increases. Hence  the vertex angles of $A_1^1$ are $\pi/4,3\pi/4,\pi/4,3\pi/4$. Let $C_3$ ($C_4$, respectively) be the great circle containing $\overline{p_1'p_2'}$ $(\overline{p_1''p_2''},$ resp.). Set $p_j^{\prime}=C_3\cap S_j^1$ ($p_j''=C_4\cap S_j^1$, resp.). Then $p_1',\ldots,p_{2k}'$ ($p_1'',\ldots,p_{2k}''$, resp.) are equally spaced on $C_3$ ($C_4$, resp.) and $\overline{p_jp_j^{\prime}}\subset T_1\cap S_j^1$ ($\overline{p_jp_j''}\subset T_2\cap S_j^1$, resp.) is an edge of the pentahedron $U_j^1$. Note that $\overline{p_j'p_j''}$ is not a geodesic in $\mathbb S^3$ but a geodesic on $T_0$: it is part of a latitude on $S_j^1$.

Let $q_1,\ldots,q_{2k}$ be the equally spaced points on $C_2$ such that $p_j'\in C_3$ is the midpoint of $\overline{p_jq_j}$, $1\leq j\leq2k$. And let $S_1^2,\ldots,S_{2k}^2$ be the great spheres containing $q_1,\ldots,q_{2k}$, respectively, and perpendicular to $C_2$. Then $T_1,\ldots,T_m$ and $S_1^2,\ldots,S_{2k}^2$ divide $D_2$ into congruent domains $\{V_j^i\}_{1\leq i\leq2m,\,1\leq j\leq2k}$.  $\{V_j^i\}$ are numbered in the same way as $\{U_j^i\}$ such that $\cup_{i=1}^{2m}V_j^i$ is a component of $D_2\sim(S_j^2\cup S_{j+1}^2)$ and $\cup_{j=1}^{2k}V_j^i,\cup_{j=1}^{2k}V_j^{i+m}$ are components of $D_2\sim(T_i\cup T_{i+1})$, and $V_j^i,V_j^{i+m}$ are symmetric about $C_2$, that is, $V_j^{i+m}=\rho_{C_2}(V_j^i)$. Again  the pentahedra $V_j^i$'s give a tessellation of $T_0$ by the parallelograms $\{B_j^i\}_{1\leq i\leq2m,1\leq j\leq2k}$. $A_j^i$ and $B_j^i$ have the same vertex angles. However, they are not congruent in $T_0$, but symmetric.

So far we know that $\{U_j^i,V_j^i\}_{1\leq i\leq2m,1\leq j\leq2k}$ forms a tessellation of $\mathbb S^3$. But we need more information than this between $U_j^i$ and $V_j^i$. Let $C_{ij}^n=T_i\cap S_j^n$, $1\leq i\leq m,1\leq j\leq 2k,n=1,2$. $C_{ij}^1$ is a great circle passing through $p_j$ and perpendicular to $C_1$, while $C_{ij}^2\ni q_j$ and $C_{ij}^2\perp C_2$. Let $q_0$ be the point on $C_2$ such that $p_2''$ is the midpoint of $\overline{p_2q_0}$ as in Figure 2. If $k$ is divisible by $m$, that is, $k=m\ell$ for some integer $\ell$, then $q_0=q_{\ell+2}$ since ${\rm dist}(q_2,q_0)=\pi/m$. Hence one can easily see that  $$C_2\cap\bigcup_{i,j}C_{ij}^1=\{q_1,\ldots,q_{2k}\}\,\,\,\,\, {\rm and}\,\,\,\,\, C_1\cap\bigcup_{i,j}C_{ij}^2=\{p_1,\ldots,p_{2k}\}.$$ It follows that \begin{equation}\label{grid}
\bigcup_{i,j}C_{ij}^1=\bigcup_{j=1}^{2k}\bigcup_{i=1}^{2m}\overline{p_jq_{j+(i-1)\ell}}=
\bigcup_{a=1}^{2k}\bigcup_{i=1}^{2m}\overline{p_{a-(i-1)\ell} q_a}=\bigcup_{i,j}C_{ij}^2.
\end{equation}
For $1\leq i\leq2m,1\leq j\leq2k$ and $b=1,2,3,4$, let $U_{jb}^i$ be the edges of $U_j^i$ perpendicular to $T_0$, and $V_{jb}^i$ those of $V_j^i$ perpendicular to $T_0$. Then one sees that \begin{equation*}
\left(\bigcup_{i,j,b}U_{jb}^i\right)\bigcup\left(\bigcup_{i,j,b}V_{jb}^i\right)=
\bigcup_{i,j}C_{ij}^1=\bigcup_{i,j}C_{ij}^2.
\end{equation*}
Hence one can conclude that $\bigcup_{i,j}C_{ij}^1$ becomes a lattice (grid) of $\mathbb S^3$ consisting of the edges of $\{U_j^i,V_j^i\}$. This observation is rather surprising, considering that the parallelograms  $A_j^i$ and $B_j^i$ are not congruent in $T_0$.

Let $\Gamma\subset\partial U_1^1$ be the piecewise geodesic Jordan curve with six ordered vertices $p_1,p_1',p_2',p_2,p_2'',p_1'',p_1$. $U_1^1$ is mean convex because it is bounded by three minimal quadrilaterals and two totally geodesic triangles. Then Jost [J] shows that $\Gamma$ spans an embedded minimal disk $H\subset U_1^1$. Define $H'=\rho _{C_3}(H)\subset V_1^{2m}=\rho_{C_3}(U_1^1)$.

Denote by $\rho_{\,\overline{pq}}$ the $180^\circ$-rotation of $\mathbb S^3$ around the great circle $\overline{pq}$. Since $H$ is bounded by six geodesic arcs $\overline{p_1p_1'},\,\overline{p_1'p_2'}, \, \overline{p_2'p_2},\,\overline{p_2p_2''},\,\overline{p_2''p_1''},\,\overline{p_1''p_1}$, $H$ can be analytically extended across the boundary by $180^\circ$-rotations.  Note that the six corresponding rotations $\rho_{\,\overline{p_1p_1'}},\rho_{\,\overline{p_1'p_2'}},\rho_{\,\overline{p_2'p_2}},\rho_{\,
\overline{p_2p_2''}},
\rho_{\,\overline{p_2''p_1''}},\rho_{\,\overline{p_1''p_1}}$ generate a finite group $G^o$ of isometries of $\mathbb S^3$. Hence one can perform those analytic extensions for all members of $G^o$ to obtain a compact minimal extension $T_{m,k}^{o}$ of $H$ without boundary. Obviously $T_{m,k}^o$ is invariant under $G^o$.

Now we claim that $T_{m,k}^{o}$ has no self intersection. Let $\bar{p}_1,\ldots,\bar{p}_{4mk}$ be the vertices of the parallelograms $A_j^i$ (such as $p_1',p_2',p_1'',$ $p_2''$), and $\bar{q}_1,\ldots,\bar{q}_{4mk}$ the vertices of $B_j^i$. Define $\rho_{\bar{p}_c}$ to be the $180^\circ$-rotation about the great circle through $\bar{p}_c$ and perpendicular to $T_0,c=1,\ldots,4mk$. Define $\rho_{\bar{q}_c}$ similarly. Extend $H$ analytically by applying $\rho_{\bar{p}_1},\ldots,\rho_{\bar{p}_{4mk}}$ to obtain $T_{m,k}^1\subset D_1$. Also extend $H'$ by applying $\rho_{\bar{q}_1},\ldots,\rho_{\bar{q}_{4mk}}$ to get $T_{m,k}^2\subset D_2$. Clearly $\partial T_{m,k}^1\subset T_0$ and $\partial T_{m,k}^2\subset T_0$. Since $T_{m,k}^1$ and $T_{m,k}^2$ are embedded, $T_{m,k}^{o}$ will be embedded if one can prove $T_{m,k}^{o}=T_{m,k}^1\cup T_{m,k}^2$. Or equivalently, $T_{m,k}^o$ will be embedded if $\rho_{C_4}(T_{m,k}^1)=T_{m,k}^2$.

Since
\begin{equation*}
T_0\cap\bigcup_{i,j}C_{i,j}^1=\{\bar{p}_1,\ldots,\bar{p}_{4mk}\}\,\,\,\,\,{\rm and}\,\,\,\,\,T_0\cap\bigcup_{i,j}C_{i,j}^2=\{\bar{q}_1,\ldots,\bar{q}_{4mk}\},
\end{equation*}
it follows from (\ref{grid}) that $$\{\bar{p}_1,\ldots,\bar{p}_{4mk}\}=\{\bar{q}_1,\ldots,\bar{q}_{4mk}\}.$$ Here we show that the divisibility of $k$ by $m$ is not sufficient for the embeddedness of $T_{m,k}^0$.

\begin{center}
\includegraphics[width=2.8in]{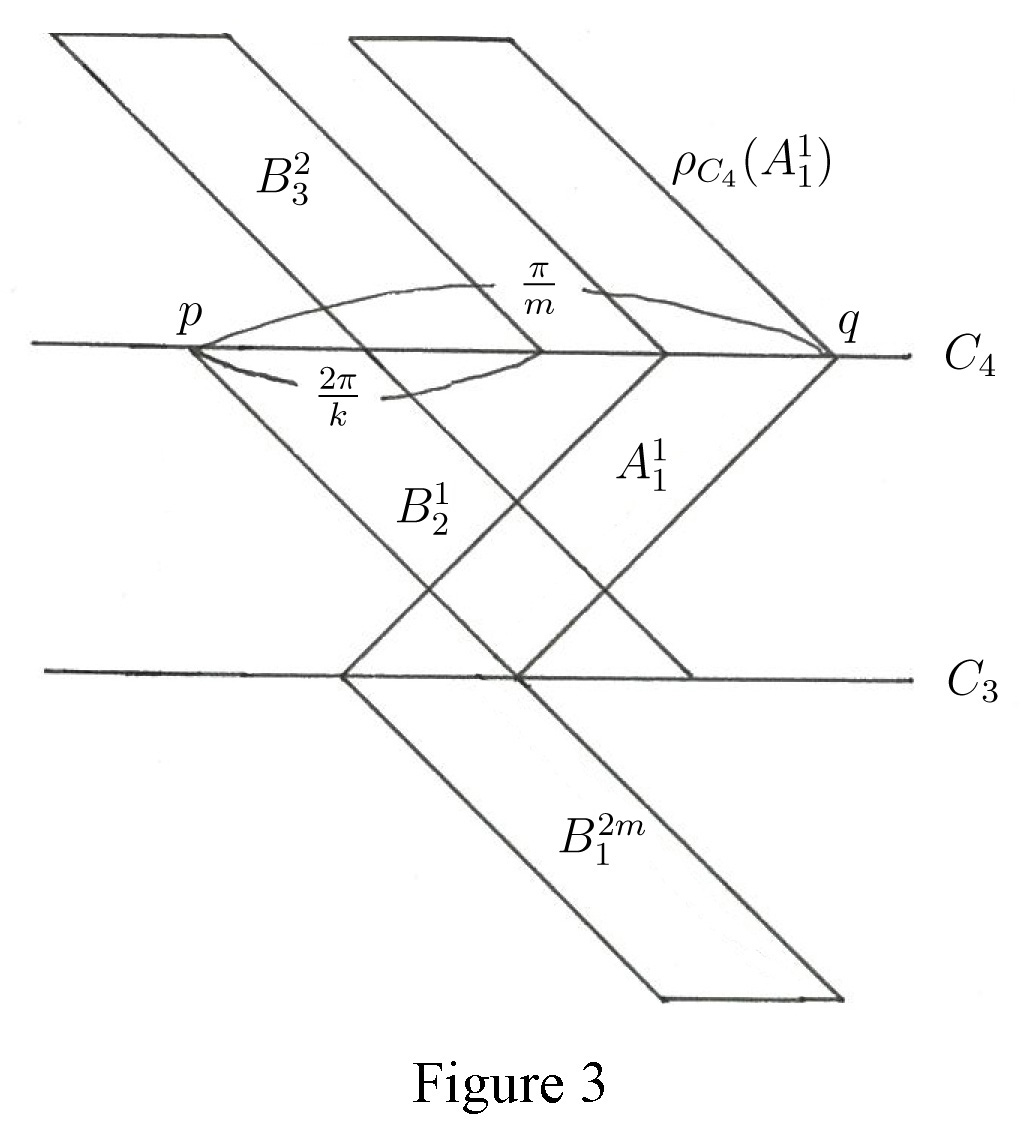}\\
\end{center}

The invariance of $T_{m,k}^1$ under the rotations $\rho_{\bar{p}_c}$ implies that $T_{m,k}^1$ occupies every other  pentahedron $U_j^i$ alternatingly. Similarly $T_{m,k}^2$ does $V_j^i$. Hence
$$T_{m,k}^1\subset\bigcup_{i+j={\rm even}}U_j^i\,\,\,\,\,\,{\rm and}\,\,\,\,\,\,T_{m,k}^2\subset\bigcup_{i+j={\rm odd}}V_j^i.$$ The length of the arc $\overline{p_2'p_2''}$ is $\frac{\pi}{\sqrt{2}m}$ in Figure 2. Since the vertex angles of $A_j^i,B_j^i$ are $\pi/4,3\pi/4,\pi/4,3\pi/4$, the length of $\overline{pq}$ is $\pi/m$ in Figure 3, $q=p_2''$. Hence if $\pi/m$ is an even multiple of $\pi/k$, that is, $k=2m\ell$ for some integer $\ell$, then one sees from Figure 3 that  $\rho_{C_4}(T_{m,k}^1)=T_{m,k}^2$. One cannot draw the same conclusion in case $k$ is an odd multiple of $m$ due to $T_{m,k}^1$'s alternating occupancy in $U_j^i$. Therefore $T_{m,k}^o=T_{m,k}^1\cup T_{m,k}^2$, and thus $T_{m,k}^{o}$ is embedded.

There are $2mk$ congruent copies of $H$ in $T_{m,k}^1$. Similarly for $T_{m,k}^2$, hence $T_{m,k}^{o}$ contains a total of $4mk$ congruent copies of $H$ when it is embedded. Now let's apply the Gauss-Bonnet theorem to $H$. Note that the external angles of $H$ are $\pi/2$ at its vertices $p_1',p_2',p_1'',p_2''$ and $(m-1)\pi/m$ at $p_1,p_2$. Hence
$$\int_HKdA+\left(4-\frac{2}{m}\right)\pi=2\pi.$$
Therefore
$$2\pi\chi(T_{m,k}^{o})=\int_{T_{m,k}^o}KdA=4mk\left(-2\pi+\frac{2\pi}{m}\right),$$
and so
$$g=1+2k(m-1).$$

For (ii) note that $\Phi_{1234}^\frac{2\pi}{k}$ maps $\overline{p_1p_1'}$ onto $\overline{p_3p_3'}$ and
$$\Phi_{1234}^\frac{2\pi}{k}(T_{m,k}^o)=T_{m,k}^o.$$
So $T_{m,k}^{o}$ is invariant under the finite cyclic group generated by $\Phi_{1234}^\frac{2\pi}{k}$.

For (iii) remember that the parallelogram $A_1^1=\square p_1'p_2'p_2''p_1''$ has vertex angles of $\pi/4,3\pi/4,\pi/4,$ $3\pi/4$. Hence the fundamental piece $H$ can have no reflection symmetry across a great sphere and neither can $T_{m,k}^{o}$. However, $T_{m,k}^{o}$ has $180^\circ$-rotation symmetries.

For (iv) note that both the minimal disk $H$ and the union $H_0$ of two flat rectangles $\square p_1p_1'p_2'p_2$ and $\square p_1p_1''p_2''p_2$ span the same Jordan curve $\Gamma$. Hence
$${\rm Area}(H)<{\rm Area}(H_0).$$ Since $$\bigcup_{\rho\in G^o}\rho(H_0)=T_1\cup\cdots \cup T_m$$
 and Area$(T_i)=2\pi^2$, the conclusion follows.
\end{pf}

\section{Even surfaces}
$T_{m,k}^{o}$ is a desingularization of $T_1\cup\cdots\cup T_m$ along $C_1\cup C_2$. But there is another way of desingularizing $T_1\cup\cdots\cup T_m$ because once $T_1\cup\cdots\cup T_m$ is desingularized along $C_1$, there are two ways of desingularization along $C_2$. The new desingularization can be done by replacing $H$ with $K$ which is obtained by freeing the edges $\overline{p_1'p_2'},\,\overline{p_1''p_2''}$ of $H$ into the curves on $A_1^1$, decreasing its area.

\begin{thm}
Let $T_1,\ldots,T_m\subset\mathbb S^3$ be the Clifford tori meeting each other along a great circle $C_1$ at an angle of $\pi/m$. Then there exists a compact minimal surface $T_{m,k}^{e}$ desingularizing $T_1\cup\cdots\cup T_m$ for each $k=2m\ell$ with integer $\ell\geq2${\rm :}

{\rm (i)} $T_{m,k}^{e}$ is embedded and has genus $1+2k(m-1)${\rm ;}

{\rm (ii)} $T_{m,k}^{e}$ is invariant under a finite group of screw motions{\rm ;}

{\rm (iii)} $T_{m,k}^{e}$ has no reflection symmetry across a great sphere{\rm ;}

{\rm (iv)} {\rm Area}$(T_{m,k}^{e})<{\rm Area}(T_{m,k}^o)$.
\end{thm}

\begin{pf}
First let's diversify the screw motion $\Phi_{ijkl}^t$. Let $C_a,C_b$ be two great circles with appropriate orientations which are polar to each other. Denote by $\rho_C^t$ the rotation of $\mathbb S^3$ around the polar great circle of $C$ by the angle $t$. Then $\rho_C^t|_C$ is the translation on $C$ by distance $t$. Now define a screw motion $\Phi_{C_a}^t$ by
$$\Phi_{C_a}^t=\rho_{C_a}^t\circ\rho_{C_b}^t.$$
We also define a screw motion with distinct speeds $\Phi_{C_a}^{t,s}$ by
$$\Phi_{C_a}^{t,s}=\rho_{C_a}^t\circ\rho_{C_b}^s.$$
Let $C_5$ be the great circle on $T_0$ which connects the midpoint of $\overline{p_1' p_1''}$ to that of $\overline{p_2'p_2''}$ and denote by $C_6$ the great circle polar to $C_5$. Put $$\varphi=\Phi_{C_5}^{-\frac{\pi}{2m},\pi}$$ and define $G^e$ to be the cyclic group generated by $\varphi$.
Clearly we have  $$\varphi(\overline{p_1''p_2''})=\overline{p_1'p_2'},\,\,\,\,{\rm and }\,\,\,\, |G^e|=4m.$$

Let $p_c\in C_5$ be the midpoint of $\overline{p_2'p_2''}$ and $p_e\in C_5$ the point closest to $p_2''$ as in Figure 5. Then $\overline{p_cp_e}$ is perpendicular to $\overline{p_ep_2''}$ and to $\overline{p_cp_2}$, and ${\rm dist}(p_c,p_e)=\pi/(4m)$. Define a Jordan curve $\Gamma$ by $\Gamma=\overline{p_2p_c}\cup\overline{p_cp_e}\cup\overline{p_ep_2''}\cup\overline{p_2''p_2}$ as in Figure 4. Obviously $\Gamma$ bounds an area minimizing disk $B_0$. $B_0$\\
\begin{center}
\includegraphics[width=2.5in]{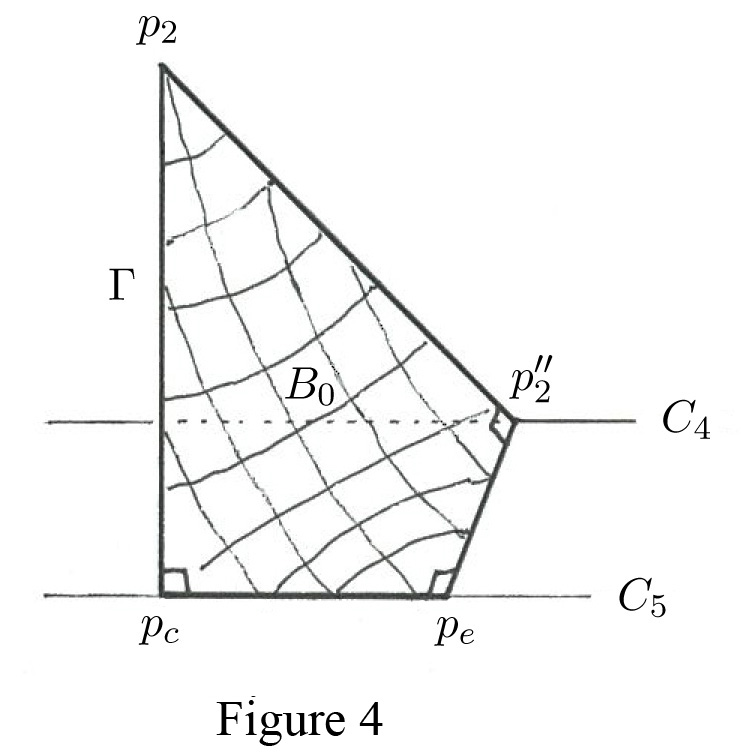}\\
\end{center}
extends analytically across $\overline{p_2p_c}$ to $B_0\cup\rho_{\,\overline{p_2p_c}}(B_0)$by the $180^\circ$-rotation $\rho_{\,\overline{p_2p_c}}$. Define
$$\Gamma_1=\bigcup_{n=1}^{4m}\varphi^n(\overline{p_2''p_2}\cup\overline{p_2p_2'}),$$
$$B_1=\bigcup_{n=1}^{4m}\varphi^n(B_0\cup \rho_{\,\overline{p_2p_c}}(B_0)),$$
$$F=\bigcup_{n=1}^{4m}\varphi^n(\square p_1p_1''p_2''p_2 \cup\,\square p_1p_1'p_2'p_2).$$
Since $\Phi_{C_1}^{-\pi/k}$ takes $p_2,p_2',p_2''$ to $p_1,p_1',p_1''$, respectively, let's also define
$$\Gamma_2=\Phi_{C_1}^{-\pi/k}(\Gamma_1)\,\,\,\,{\rm and}\,\,\,\,B_2=\Phi_{C_1}^{-\pi/k}(B_1).$$
Then $\Gamma_1$, $\Gamma_2$ are  helical Jordan curves consisting of $4m$ geodesic arcs and winding around $C_5$, and $B_1$, $B_2$ are embedded minimal annuli spanning $\Gamma_1\cup C_5$, $\Gamma_2\cup C_5$, respectively. Let $W$ be the domain bounded by $F\cup B_1\cup B_2$. Obviously $\Gamma_1,\Gamma_2,B_1,B_2,F$ and $W$ are all invariant under the cyclic group $G^e$.

We now claim that $W$ is mean convex and that there exists an embedded minimal annulus $K^e$ in $W$ spanning $\Gamma_1\cup\Gamma_2$.  For the mean convexity of $W$ it suffices to show that $ B_1$ and $B_2$ make an angle $\leq\pi$ along their intersection, $C_5$. Since $B_1$ and $B_2$ are invariant under $G^e$, we have only to prove this angle condition on an arc ($=\overline{p_ap_d}$) of length $\pi/(2m)$ in $C_5$.
\begin{center}
\includegraphics[width=3.5in]{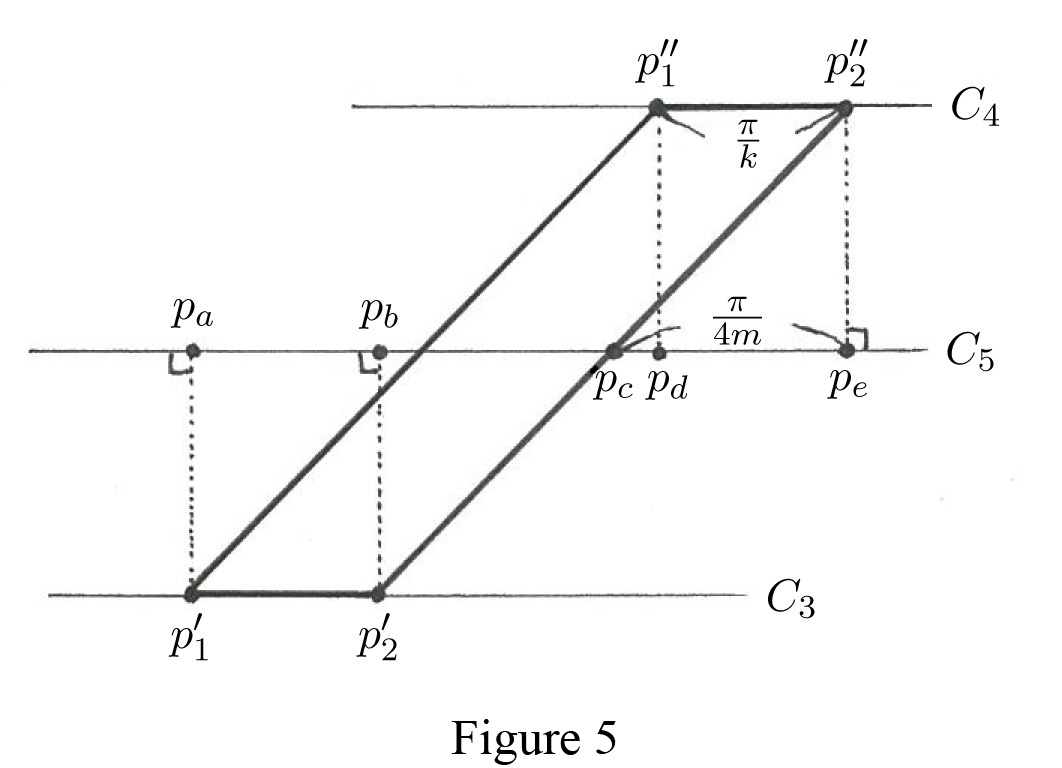}\\
\end{center}

Let $p_a,p_b,p_d$ be the points of $C_5$ closest to $p_1',p_2',p_1''$, respectively (see Figure 5). Then ${\rm dist}(p_a,p_d)={\rm dist}(p_b,p_e)=\pi/(2m)$. It is here that we need the hypothesis $k\geq4m$, i.e., $\ell\geq2$. Then
\begin{equation}\label{hypo}
{\rm dist}(p_a,p_b)\leq{\rm dist}(p_b,p_d).
\end{equation}
On $\overline{p_bp_d}$ both $B_1$ and $B_2$ are on the same side of $T_0$ because $\Phi_{C_1}^{-\pi/k}(T_0)=T_0$. Hence they make an angle $\leq\pi$ along $\overline{p_bp_d}$. On the other hand, note that along $\overline{p_cp_e}$ $B_0$ makes an acute angle with the component $T_0^4$ of $T_0\sim(C_4\cup C_5)$ containing $\overline{p_ep_2''}$ (see Figure 4). Hence along $\overline{p_ap_b}$ (\ref{hypo}) implies that  $B_2$ makes an acute angle with $T_0^3$, where $T_0^3$ is the component of $T_0\sim (C_3\cup C_5)$ containing $\overline{p_2'p_c}$. Moreover $T_0^3$ makes an acute angle with $\varphi(B_0)\subset B_1$ along $\overline{p_ap_b}$. Therefore $B_1$ and $B_2$ make an angle $\leq\pi$ along $\overline{p_ap_b}$. So $W$ is locally mean convex along $\overline{p_ap_d}$, and it follows that $W$ is mean convex. Let $U$ be the component of $W\sim T_0$ such that $\bar{U}\supset\overline{p_1p_2}$. Then $U$ is also mean convex.

Denote by $\mathcal{A}$ the set of all curves $\alpha\subset \partial U \cap T_0^4$ from $p_1''$ to $p_2''$ with no self intersection. For $\alpha\in \mathcal{A}$, let $\Gamma_{\alpha}$ be the Jordan curve $\alpha\cup\,\overline{p_2''p_2}\cup\,\overline{p_2p_2'}\cup\,\varphi(\alpha)\cup\,
\overline{p_1'p_1}\cup\,\overline{p_1p_1''}$. The mean-convexity of $U$ guarantees the existence of an embedded minimal disk in $U$ spanning $\Gamma_\alpha$. Let $\mathcal{K}$ be the family of all such embedded minimal disks in $U$ spanning $\Gamma_{\alpha}$ for all $\alpha\in \mathcal{A}$.

We now show that there exists an area minimizer $K$ in $\mathcal{K}$:
$${\rm Area}(K)=\inf_{\hat{K}\in\mathcal{K}}{\rm Area}(\hat{K}).$$
Let $\{\hat{K_i}\}$ be a minimizing sequence in $\mathcal{K}$ with
$$\lim_{i\rightarrow\infty}{\rm Area}(\hat{K_i})=\inf_{\hat{K}\in\mathcal{K}}{\rm Area}(\hat{K}).$$
The limit $K$ of $\{\hat{K}_i\}$ as an area minimizing current obviously exists in $U$. And it is easy by \cite{J} to see that $K$ is a smooth embedded minimal disk.

Denote again by $\alpha$ the part of $\partial K\,\cap\, T_0^4$ connecting $p_1''$ to $p_2''$. Then $\alpha\,\cup\,\varphi(\alpha)=\partial K\cap T_0$, and one can see that $K$  analytically extends to $K\cup\,\varphi(K)$ across $\varphi(\alpha)$. This is because if $K$  makes an angle $\neq\pi$ with $\varphi(K)$ along $\varphi(\alpha)$, then $K\cup\varphi(K)$ can be perturbed along $\varphi(\alpha)$ decreasing its area, which contradicts the assumption that $K$ is a minimizer. Similarly $K$ should extend analytically to $K\cup(\varphi)^{-1}(K)$ across $\alpha$. Therefore one can extend $K$ analytically to a smooth minimal annulus
$$K^e:=\bigcup_{n=1}^{4m}\varphi^n(K)~~~~~\subset W$$
which spans $\Gamma_1\cup\Gamma_2$. Note that $$\Gamma_1\cup\Gamma_2\subset\bigcup_{i,j}C_{ij}^1.$$ Therefore $K^e$ can be indefinitely extended across $\Gamma_1\cup\Gamma_2$ by $180^\circ$-rotations $\rho_{\bar{p}_1},\ldots,\rho_{\bar{p}_{4mk}}$ to produce a complete minimal surface $T_{m,k}^{e}$. Since the group generated by $\rho_{\bar{p}_1},\ldots,\rho_{\bar{p}_{4mk}}$ is finite, one can conclude that $T_{m,k}^{e}$ is compact.

Let's show that $T_{m,k}^{e}$ is embedded. Clearly $\varphi(U_1^1)=V_1^1$ and $\varphi(A_1^1)=B_1^1$ as in Figure 6. Suppose $k$ is divisible by $2m$, i.e., $k=2m\ell$. Then $$\varphi^2(U_1^1)=U_{2k-2\ell+1}^1\subset \bigcup_{i+j={\rm even}}U_j^i\,\,\,\,\,{\rm and}\,\,\,\,\,\varphi^2(V_1^1)=V_{2k-2\ell+1}^1\subset \bigcup_{i+j={\rm even}}V_j^i.$$ Since $K^e$ is invariant under $\varphi ^n$, $1\leq n\leq4m$, one sees that
$$K^e\cap D_1\subset \bigcup_{i+j={\rm even}}U_j^i\,\,\,\,\,{\rm and}\,\,\,\,\,K^e\cap D_2\subset \bigcup_{i+j={\rm even}}V_j^i.$$ Hence the embeddedness of $T_{m,k}^e$ follows from the fact that for $c=1,\ldots,4mk$
$$\rho_{\bar{p}_c}\left(\bigcup_{i+j={\rm even}}U_j^i\right)=\bigcup_{i+j={\rm even}}U_j^i\,\,\,\,\,{\rm and}\,\,\,\,\,\rho_{\bar{p}_c}\left(\bigcup_{i+j={\rm even}}V_j^i\right)=\bigcup_{i+j={\rm even}}V_j^i.$$

There are $2mk$ congruent copies of $K$ in $T_{m,k}^{e}\cap D_1$, and there are the same number of copies of $\varphi(K)$ in $T_{m,k}^{e}\cap D_2$. It should be remarked that the sum of the geodesic curvatures of $\alpha$ and $\varphi(\alpha)$ at $p\in\alpha$ and $\varphi(p)\in\varphi(\alpha)$, respectively, vanishes. And the external angles of $K$ at $p_1',p_2',p_1'',p_2''$ are $\pi/2$ and $(m-1)\pi/m$ at $p_1,p_2$. So by the same argument as in Theorem 1 we see that $$g=1+2k(m-1).$$

\begin{center}
\includegraphics[width=3.5in]{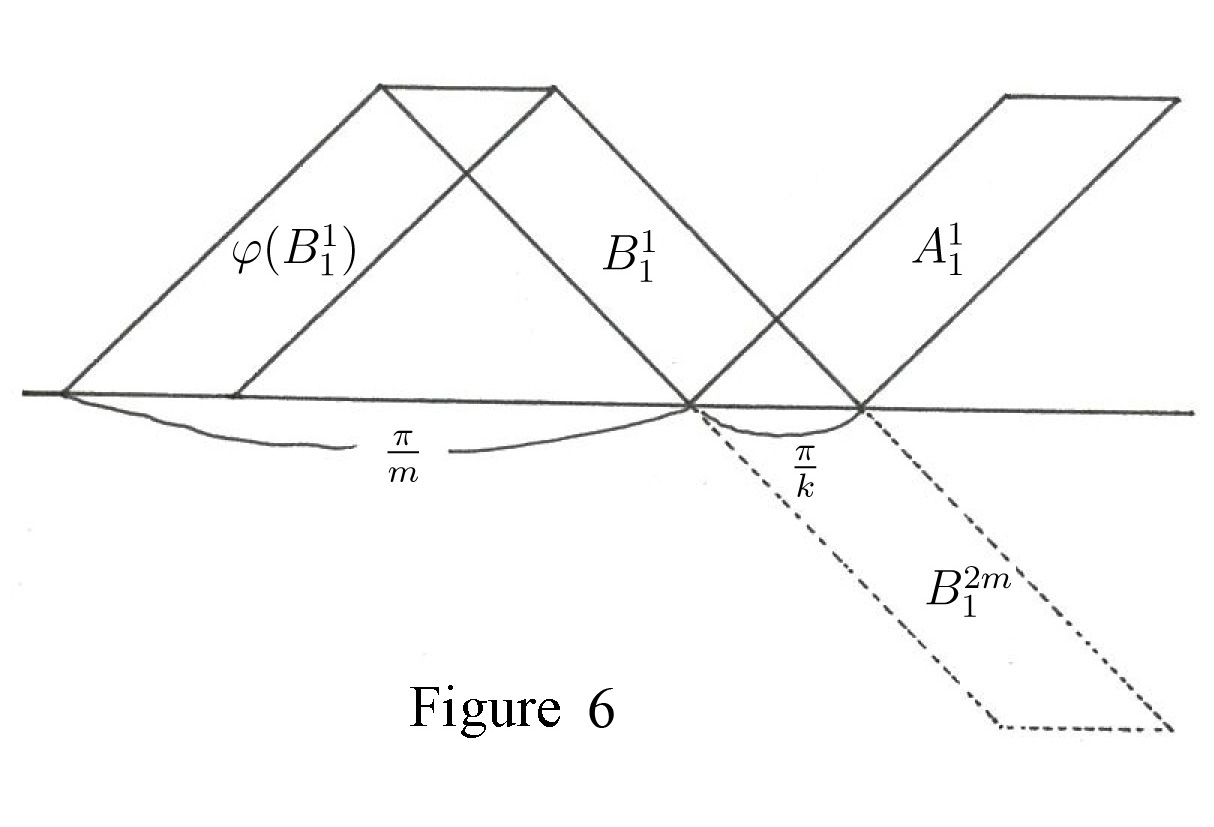}\\
\end{center}

(ii) and (iii) follow from the same arguments as in Theorem 1. For (iv) just note that $H\in \mathcal{K}$.
\end{pf}

\end{document}